\documentclass[12pt]{amsart}
\usepackage{amsmath,amsfonts,amssymb,amsthm}
\usepackage{graphicx,color}
\usepackage{hyperref}

\DeclareMathRadical{\sqrtsign}{symbols}{"70}{largesymbols}{"70}

\newlength{\figboxwidth}             
\newcommand{\infinity}{\infty}

\newtheorem{theorem}{Theorem}[section]
\newtheorem*{mydef}{Theorem 1}

\newtheorem*{mydef4}{Theorem 2}

\newtheorem{lemma}[theorem]{Lemma}

\theoremstyle{definition}
\newtheorem{definition}[theorem]{Definition}

\theoremstyle{remark}
\newtheorem{remark}[theorem]{Remark}
\theoremstyle{proposition}
\newtheorem{proposition}{Proposition}
\theoremstyle{Corollary}
\newtheorem{corollary}{Corollary}

\numberwithin{equation}{section}

\newcommand{\abs}[1]{\lvert#1\rvert}

	\def\Diff{\operatorname{Diff}}
	\def\inte{\operatorname{int}}

\begin{document}
	
\title[Continuity in generic dynamical spectra]
{Continuity of fractal dimensions in conservative generic Markov and Lagrange dynamical spectra}

\author{Davi Lima}
\address{Davi Lima: Instituto de Matem\'atica, Universidade Federal de Alagoas, Av Lourival Melo Mota s/n, 57072-970, Macei\'o-Brazil}
\email{davi.santos@im.ufal.br}
\thanks{The first author was partially supported by CNPq and FAPEAL}

\author{Carlos Gustavo Moreira}
\address{Carlos Gustavo Moreira: SUSTech International Center for Mathematics, Shenzhen, Guangdong, People’s Republic of China and IMPA, Estrada Dona Castorina 110, 22460-320, Rio de Janeiro, Brazil
}
\email{gugu@impa.br}
\thanks{The second author was partially supported by CNPq and FAPERJ}

\author{Christian Villamil}
\address{Christian Villamil: IMPA, Estrada Dona Castorina 110, 22460-320, Rio de Janeiro, Brazil}
\email{ccsilvav@impa.br}
\thanks{The third author was partially supported by CNPq}

\date{\today}


\begin{abstract}
Let $\varphi_0$ be a smooth conservative diffeomorphism of a compact surface $S$ and let $\Lambda_0$ be a transitive horseshoe of $\varphi_0$. Given a smooth real function $f$ defined in $S$ and a small smooth conservative perturbation $\varphi$ of $\varphi_0$, let $L_{\varphi, f}$ and $M_{\varphi, f}$ be respectively the Lagrange and Markov spectra associated to the  hyperbolic continuation $\Lambda(\varphi)$ of the horseshoe $\Lambda_0$ and $f$. We show that for generic choices of $\varphi$ and $f$, the Hausdorff dimension of the sets $L_{\varphi, f}\cap (-\infty, t)$ and $M_{\varphi, f}\cap (-\infty, t)$ are equal and determine a continuous function as $t\in \mathbb{R}$ varies; generalizing then the Cerqueira-Matheus-Moreira theorem to horseshoes with arbitrary Hausdorff dimension.  
\end{abstract}

\maketitle

\tableofcontents

\section{Introduction}

\subsection{Classical spectra}

The classical Lagrange and Markov spectra are closed subsets of the real line related to Diophantine approximations. They arise naturally in the study of rational approximations of irrational numbers and of indefinite binary quadratic forms, respectively. More precisely, given an irrational number $\alpha$, let 
$$k(\alpha):=\limsup_{\substack{p, q\to\infty \\ p, q\in\mathbb{N}}}(q|q\alpha-p|)^{-1}$$
be its best constant of Diophantine approximation. The set 
$$L:=\{k(\alpha):\alpha\in\mathbb{R}-\mathbb{Q}, k(\alpha)<\infty\}$$ 
consisting of all finite best constants of Diophantine approximations is the so-called \emph{Lagrange spectrum}. 

Similarly, the \emph{Markov spectrum} 
 $$M:=\left\{\left(\inf\limits_{(x,y)\in\mathbb{Z}^2-\{(0,0)\}} |q(x,y)|\right)^{-1} < \infty: q(x,y)=ax^2+bxy+cy^2, b^2-4ac=1\right\} $$
consists of the reciprocal of the minimal values over non-trivial integer vectors $(x,y)\in\mathbb{Z}^2-\{(0,0)\}$ of indefinite binary quadratic forms $q(x,y)$ with unit discriminant. 

Hurwitz showed that the minimum of $L$ is $\sqrt{5}.$ Markov improved the result by showing that 
$$L\cap (-\infty, 3)=M\cap (-\infty, 3)=\{k_1=\sqrt{5}<k_2=2\sqrt{2}<k_3=\frac{\sqrt{221}}{5}<...\},$$
where $k^2_n\in \mathbb{Q}$ for every $n\in \mathbb{N}$ and $k_n\to 3$ when $n\to \infty$.

Hall showed that $L$ contains a half-line $[c,+\infty)$ for some $c>4$, and Freiman determined the biggest such a half-line.

Moreira in \cite{M3} proved several results on the geometry of the Markov and Lagrange spectra, for example, that the map $d:\mathbb{R} \rightarrow [0,1]$, given by
$$
d(t)=HD(L\cap(-\infty,t))= HD(M\cap(-\infty,t))
$$
is continuous, surjective and such that $d(3)=0$ and $d(\sqrt{12})=1$. Moreover, that
$$d(t)=\min \{1,2D(t)\}$$
where $D(t)=HD(k^{-1}(-\infty,t))=HD(k^{-1}(-\infty,t])$ is also a continuous function from $\mathbb{R}$ to $[0,1).$ Even more, he proved the limit
$$\lim_{t\rightarrow \infty}HD(k^{-1}(t))=1.$$

In this work, we will be mainly interested in that kind of result. The reader can find more information about the structure of these sets in the classical book \cite{CF} of Cusick and Flahive, but for our purposes, it is worth to point out here that the Lagrange and Markov spectra have the following \emph{dynamical} interpretation in terms of the continued fraction algorithm: 

Denote by $[a_0,a_1,\dots]$ the continued fraction $a_0+\frac{1}{a_1+\frac{1}{\ddots}}$. Let $\Sigma=\mathbb{N}^{\mathbb{Z}}$ the space of bi-infinite sequences of positive integers, $\sigma:\Sigma\to\Sigma$ be the left-shift map $\sigma((a_n)_{n\in\mathbb{Z}}) = (a_{n+1})_{n\in\mathbb{Z}}$, and let $f:\Sigma\to\mathbb{R}$ be the function
$$f((a_n)_{n\in\mathbb{Z}}) = [a_0, a_1,\dots] + [0, a_{-1}, a_{-2},\dots].$$
Then, 
$$L=\left\{\limsup_{n\to\infty}f(\sigma^n(\underline{\theta})):\underline{\theta}\in\Sigma\right\} \quad \textrm{and} \quad M= \left\{\sup_{n\to\infty}f(\sigma^n(\underline{\theta})):\underline{\theta}\in\Sigma\right\}.$$
It follows from these characterizations that $M$ and $L$ are closed subsets of $\mathbb R$ and that $L\subset M$.

In the sequel, we consider the natural generalization of this dynamical version of the classical Lagrange and Markov spectra in the context of horseshoes\footnote{i.e., a non-empty compact invariant hyperbolic set of saddle type which is transitive, locally maximal, and not reduced to a periodic orbit (cf. \cite{PT93} for more details).} of smooth diffeomorphisms of compact surfaces. In this setting, our main result (cf. Theorem \ref{main2} below) will be a generalization of the results of \cite{CMM16} on the continuity of Hausdorff dimension across Lagrange and Markov dynamical spectra.

\subsection{Dynamical spectra}

Let $\varphi:S\rightarrow S$ be a diffeomorphism of a $C^{\infty}$ compact surface $S$ with a transitive horseshoe $\Lambda$ and let $f:S\rightarrow \mathbb{R}$ be a differentiable function. Following the above characterization of the classical spectra, we define the maps
\begin{eqnarray*}
   \ell_{\varphi,f}: S &\rightarrow& \mathbb{R} \\
   x &\mapsto& \ell_{\varphi,f}(x)=\limsup_{n\to \infty}f(\varphi^n(x))
\end{eqnarray*}
\begin{eqnarray*}
   m_{\varphi,f}:S &\rightarrow& \mathbb{R} \\
   x &\mapsto& m_{\varphi,f}(x)=\sup_{n\to \infty}f(\varphi^n(x))
\end{eqnarray*} and

and call $\ell_{\varphi,f}(x)$ \textbf{the Lagrange value} of $x$ associated to $f$ and $\varphi$ and also $m_{\varphi,f}(x)$ \textbf{the Markov value} of $x$ associated to $f$ and $\varphi$. The sets
$$L_{\varphi,f}=\ell_{\varphi,f}(\Lambda)=\{\ell_{\varphi,f}(x):x\in \Lambda\}$$
and
$$M_{\varphi,f}=m_{\varphi,f}(\Lambda)=\{m_{\varphi,f}(x):x\in \Lambda\}$$
are called \textbf{Lagrange Spectrum} of $(\varphi,f,\Lambda)$ and \textbf{Markov Spectrum} of $(\varphi,f, \Lambda)$.
An elementary compactness argument (cf. Remark in Section 3 of \cite{MR2}) shows that  $\{\ell_{\varphi, f}(x):x\in X\}\subset\{m_{\varphi, f}(x): x\in X\}\subset f(X)$
whenever $X\subset S$ is a compact $\varphi$-invariant subset.

If $HD(X)$ denotes the Hausdorff dimension of $X$, in this paper we are interested in the study of the two real functions
\begin{equation}\label{f1}
L(t)=L(\varphi,f)(t)=HD(L_{\varphi,f}\cap (-\infty,t))
\end{equation}
and
\begin{equation}
M(t)=M(\varphi,f)(t)=HD(M_{\varphi,f}\cap (-\infty,t)).
\end{equation}
For this reason, we will also study the fractal geometry (Hausdorff dimension) of 
$$\Lambda_t:=\bigcap\limits_{n\in\mathbb{Z}}\varphi^{-n}(\{y\in\Lambda: f(y)\leq t\}) = \{x\in\Lambda: m_{\varphi, f}(x)=\sup\limits_{n\in\mathbb{Z}}f(\varphi^n(x))\leq t\}$$ 
for $t\in\mathbb{R}$. 

Fix a Markov partition $\{R_a\}_{a\in \mathcal{A}}$ of $\Lambda$ with sufficiently small diameter consisting of rectangles $R_a \sim I_a^s \times I_a^u$ delimited by compact pieces $I_a^s$, $I_a^u$, of stable and unstable manifolds of certain points of $\Lambda$, see \cite{PT93} theorem 2, page 172. Then we can define the subset $\mathcal{B}\subset \mathcal{A}^2$ of admissible transitions as the subset of pairs $(a,b)$ such that $\varphi(R_a)\cap R_b\neq \emptyset$. 
In this way, there is a homeomorphism $\Pi:\Lambda\rightarrow \Sigma$ such that 
	\begin{equation}\label{conjugated}
	\Pi(\varphi(x))=\sigma(\Pi(x))
	\end{equation}
	where $\Sigma=\Sigma_{\mathcal{B}}$ is the Markov shift of finite type associated to $\mathcal{B}$ and $\sigma$ is the left-shift map as before. 
	
	Next, we recall that the stable and unstable manifolds of $\Lambda$ can be extended to locally invariant $C^{1+\varepsilon}$ foliations in a neighborhood of $\Lambda$ for some $\varepsilon>0$. Therefore, we can use these foliations to define projections $\pi^u_a: R_a\rightarrow I^s_a \times \{i^u_a\}$ and $\pi^s_a: R_a\rightarrow \{i^s_a\}\times I^u_a$ of the rectangles into the connected components $I^s_a \times \{i^u_a\}$ and $\{i^s_a\}\times I^u_a$ of the stable and unstable boundaries of $R_a$, where $i^u_a\in \partial I^u_a$ and $i^s_a\in \partial I^s_a$ are fixed arbitrarily. Using these projections, we have the unstable and stable Cantor sets
	$$K^u=\bigcup_{a\in \mathcal{A}}\pi^s_a(\Lambda\cap R_a) \ \mbox{and} \ K^s=\bigcup_{a\in \mathcal{A}}\pi^u_a(\Lambda\cap R_a).$$
\begin{figure}[h]
\centering
\includegraphics[width=0.8 \textwidth]{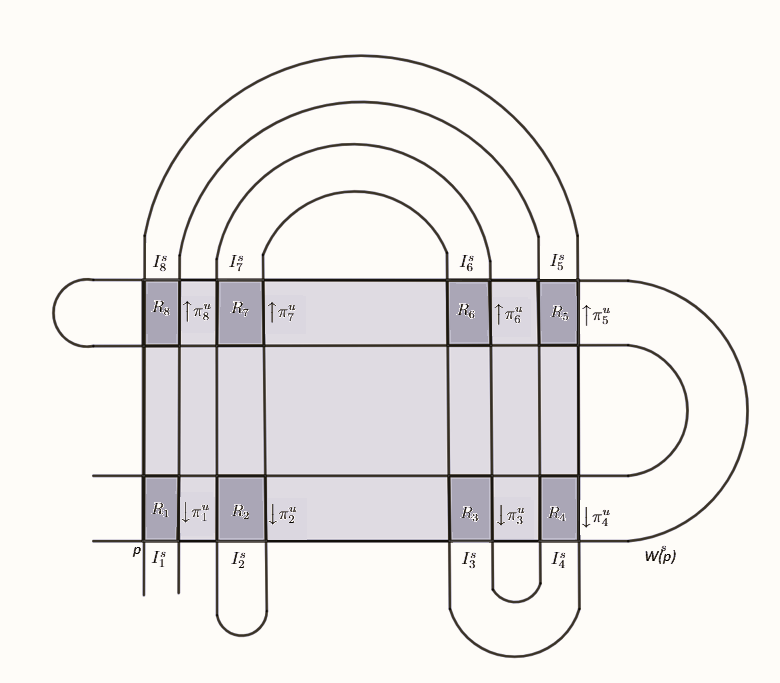}
\caption{Markov partition and projections.}
\end{figure}

\begin{remark}\label{remarkA}
In what follows, formally speaking, in order to show that $K^u$ and $K^s$ are regular Cantor sets (which, by definition, should be defined by a mixing expanding dynamics) we should suppose the horseshoe to be mixing instead of only transitive, but we don't loose generality since, if $\Lambda$ is transitive but not mixing then, by spectral decomposition, there is an integer $k>1$ and disjoint compact sets $\Lambda^j\subset\Lambda, 0\le j\le k-1$ such that $\Lambda=\bigcup_{j=0}^{k-1} \Lambda^j$, $\varphi(\Lambda^j)=\Lambda^{j+1 \pmod k}, 0\le j\le k-1$ and $\Lambda^j$ is a mixing horseshoe for the diffeomorphism $\varphi^k$, which is the first return map of the dynamics of $\varphi$ to $\Lambda^j$. So, in the following discussion, we may replace $\varphi$ by $\varphi^k$ in a neighbourhood of each $\Lambda_j$, and we may assume that each rectangle $R_a$ of the Markov partition is contained in one of these neighbourhoods. 
\end{remark}

In fact $K^u$ and $K^s$ are $C^{1+\varepsilon}$ dynamically defined Cantor sets. We define $g_s$ and $g_u$ in the following way: If $y\in R_{a_1}\cap \varphi(R_{a_0})$ we put
$$g_s(\pi^u_{a_1}(y))=\pi^u_{a_0}(\varphi^{-1}(y))$$
and if $z\in R_{a_0}\cap \varphi^{-1}(R_{a_1})$ we put
$$g_u(\pi^s_{a_0}(z))=\pi^s_{a_1}(\varphi(z)).$$

We have that $g_s$ and $g_u$ are $C^{1+\varepsilon}$ expanding maps of type $\Sigma_{\mathcal{B}}$ defining $K^s$ and $K^u$ in the sense that

\begin{enumerate}
	\item[(i)] The domains of $g_s$ and $g_u$ are disjoint unions $$\bigsqcup_{(a_0,a_1)\in \mathcal{B}} I^s(a_1,a_0) \ \mbox{and} \bigsqcup_{(a_0,a_1)\in \mathcal{B}}I^u(a_0,a_1),$$ where $I^s(a_1,a_0)$, resp. $I^u(a_0,a_1)$, are compact subintervals of $I^s_{a_1}$, resp. $I^u_{a_0}$;
	
	\item[(ii)] For each $(a_0,a_1)\in \mathcal{B},$ the restrictions $g_s|_{I^s(a_1,a_0)}$ and $g_u|_{|I^u(a_0,a_1)}$ are $C^{1+\epsilon}$ diffeomorphisms onto $I^s_{a_0}$ and $I^u_{a_1}$ with $|Dg_s(t)|,|Dg_u(t)|>1$, for all $t\in I^s(a_1,a_0)$, $t\in I^u(a_0,a_1)$ (for appropriate choices of the parametrization of $I^s_a$ and $I^u_a$);  
	
	\item[(iii)] $K^s$ and $K^u$ satisfies 
	$$K^s=\bigcap_{n\ge 0}g_s^{-n}\left(\bigsqcup_{(a_0,a_1)\in \mathcal{B}} I^s(a_1,a_0)\right) \ \ \ K^u=\bigcap_{n\ge 0}g_u^{-n}\left(\bigsqcup_{(a_0,a_1)\in \mathcal{B}}I^u(a_0,a_1)\right).$$
	
\end{enumerate}

The stable and unstable Cantor sets, $K^s$ and $K^u$, respectively, are closely related to the fractal geometry of the horseshoe $\Lambda$; for instance, it is well-known that, 
\begin{equation}\label{eq.1}
HD(\Lambda)=HD(K^s)+HD(K^u),
\end{equation}
see \cite{MMan} theorem 2 or \cite{PT93} proposition 4, page 75. 

Following the above construction we will study the subsets $\Lambda_t$ introduced above through its projections on the stable and unstable Cantor sets of $\Lambda$:
$$K^u_t=\bigcup_{a\in \mathcal{A}} \pi^s_a(\Lambda_t\cap R_a) \ \mbox{and} \ K^s_t=\bigcup_{a\in \mathcal{A}}\pi^u_a(\Lambda_t\cap R_a).$$

Also, we define in the context of transitive horseshoes $\Lambda$ with $HD(\Lambda)>1$ the Markov transition parameter as
$$a=a(\varphi,f)=\sup \{t\in \mathbb{R}:HD(\Lambda_t)<1\}.$$
In \cite{LM} is proved that for typical choices of the diffeomorphism $\varphi$ and the smooth real map $f$, the Markov parameter is characterized by the conditions
$$Leb(M_{\varphi,f}\cap (-\infty,a-\delta))=0$$
but
$$Int(M_{\varphi,f}\cap (-\infty,a+\delta))\neq \emptyset,$$
for all $\delta>0.$\footnote{\ here $Leb(\cdot)$ denotes the usual Lebesgue measure and $Int(\cdot)$ the interior of the set.}

The Lagrange parameter $\tilde{a}=\tilde{a}(\varphi,f)$ is defined in such a way that a similar result is true if we replace $M_{\varphi,f}$ by $L_{\varphi,f}$ and $a$ by $\tilde{a}$ in the last conditions. Since $L_{\varphi,f} \subset M_{\varphi,f}$, we have always $a(\varphi,f)\leq \tilde{a}(\varphi,f)$.

\subsection{Acknowledgements}
 The first author is partially supported by CNPq grant - Alagoas Din\^amica 409198/2021-8 and FAPEAL - grant E60030.0000002330/2022.

\section{Statement of the results}

The aim of this work is to extend the main theorem in \cite{CMM16}, removing the hypothesis that $HD(\Lambda) < 1$. Using the notations of the previous subsection, our results are the following.

\begin{mydef}\label{main11}
Let $\varphi\in \Diff^{2}(S)$ with a transitive horseshoe $\Lambda$. For every $r\ge 2$ there exists a $C^r$-open and dense set $\mathcal{R}_{\varphi,\Lambda}$ such that for any function $f\in \mathcal{R}_{\varphi,\Lambda}$ the functions
	$$t\mapsto d_u(t):= HD(K^u_t) \ \mbox{and} \ t\mapsto d_s(t):=HD(K^s_t)$$
are continuous.
\end{mydef}

\begin{remark}
Our proof of theorem 1 shows that $d_u$ and $d_s$ coincide with the box-counting dimensions of $K^s_t$ and $K^u_t$ respectively.
\end{remark}

We write $\Diff^{2}_{\omega}(S)$ for the set of conservative diffeomorphisms of $S$ with respect to a volume form $\omega$. Then we have the
\begin{mydef4}\label{main2}

Let $\varphi_0\in \Diff^{2}_{\omega}(S)$ with a transitive horseshoe $\Lambda_0$ and $\mathcal{U}$ a $C^{2}$-sufficiently small neighbourhood of $\varphi_0$ in $\Diff^{2}_{\omega}(S)$ such that $\Lambda_0$ admits a continuation $\Lambda (= \Lambda(\varphi))$ for every $\varphi \in \mathcal{U}$. There exists a residual set $\tilde{\mathcal{U}}\subset \mathcal{U}$ such that for every $\varphi\in \tilde{\mathcal{U}}$ and $r\ge 2$ there exists a $C^r$-residual set $\tilde{\mathcal{R}}_{\varphi,\Lambda}\subset C^r(S,\mathbb{R})$ such that for any $f\in \tilde{\mathcal{R}}_{\varphi,\Lambda}$ the functions:
	$$t\mapsto d_u(t):= HD(K^u_t) \ \mbox{and} \ t\mapsto d_s(t):=HD(K^s_t)$$
are continuous and in fact, they are equal with
$$HD(\Lambda_t)=d_u(t)+d_s(t)=2d_u(t)$$ and 
	$$\min\{1,HD(\Lambda_t)\}=L(t)=M(t).$$
\end{mydef4}

Finally, in theorem D of \cite{LM} is shown in the conservative case, that generically we have the equality $a=\tilde{a}$ where $a=a(\varphi,f)$ and $\tilde{a}=\tilde{a}(\varphi,f)$ are as in the last section. However, there is a mistake in the proof of the last statement of that theorem; more specifically, in the proof of the affirmation  
$$HD(M_{\varphi,f}\cap (-\infty,a))=HD(L_{\varphi,f}\cap (-\infty,a))=1.$$
Nevertheless, working in the setting of theorem \ref{main2} we have
\begin{eqnarray*}
       L(a)=M(a)=\min \{1, HD(\Lambda_a) \}&=&\lim \limits_{t \rightarrow a^-}\min \{1, HD(\Lambda_t) \} \\ &=&\lim \limits_{t \rightarrow a^-}HD(\Lambda_t)\\ &=&HD(\Lambda_a)=1 
\end{eqnarray*}
 then, intersecting the residual sets of the theorem D with the residual sets that we obtained here, we get a correct proof of the

\begin{corollary} [Theorem D of \cite{LM}]

Let $\varphi_0\in \Diff^{2}_{\omega}(S)$ with a transitive horseshoe $\Lambda_0$ with $HD(\Lambda_0)>1$ and $\mathcal{V}$ a $C^{2}$-sufficiently small neighbourhood of $\varphi_0$ in $\Diff^{2}_{\omega}(S)$ such that $\Lambda_0$ admits a continuation $\Lambda$ for every $\varphi \in \mathcal{V}$. Then, there exists a residual set $\mathcal{V}^*\subset \mathcal{V}$ such that for every $\varphi\in \mathcal{V}^*$ and $r\ge 2$ there exists a $C^r$-residual set $\mathcal{P}_{\varphi,\Lambda}\subset C^r(M,\mathbb{R})$ such that for any $f\in \mathcal{P}_{\varphi,\Lambda}$:
   $$Leb(M_{\varphi,f}\cap (-\infty,a-\delta))=0=Leb(M_{\varphi,f}\cap (-\infty,a-\delta))$$
   but 
   $$Int(M_{\varphi,f}\cap (-\infty,a+\delta))\neq \emptyset \neq Int(L_{\varphi,f}\cap (-\infty,a+\delta))$$
for all $\delta>0$.
Moreover, one has
$$HD(M_{\varphi,f}\cap (-\infty,a))=HD(L_{\varphi,f}\cap (-\infty,a))=1.$$

\end{corollary}
 
\section{Preliminary results }\label{s1}

Fix a Markov partition $\mathcal{P}=\{R_a\}_{a\in \mathcal{A}}$; recall that the geometrical description of $\Lambda$ in terms of the Markov partition $\mathcal{P}$ has a combinatorial counterpart in terms of the Markov shift $\Sigma_{\mathcal{B}}\subset \mathcal{A}^{\mathbb{Z}}$. And we can use $\Pi$ (see \ref{conjugated}) to transfer the function $f$ from $\Lambda$ to a function (still denoted $f$) on $\Sigma_{\mathcal{B}}$. In this setting, $\Pi(\Lambda_t)=\Sigma_t$ where 
$$\Sigma_t=\{\theta\in\Sigma_{\mathcal{B}}: \sup\limits_{n\in\mathbb{Z}} f(\sigma^n(\theta))\leq t\}.$$

Given an admissible finite sequence $\theta=(a_1,...,a_n)\in \mathcal{A}^n$ (i.e., $(a_i,a_{i+1})\in \mathcal{B}$) for all $1\le i<n$, we define
	$$I^u(\theta)=\{x\in K^u: g_u^i(x)\in I^u(a_i,a_{i+1}), i=1,2,...,n-1\}$$
and
	$$I^s(\theta^t)=\{y\in K^s: g_s^i(y)\in I^s(a_{i},a_{i-1}), i=2,...,n\}$$
 where $\theta^t=(a_n,a_{n-1},...,a_2,a_1)$. In a similar way, let $\alpha=(a_{s_1},a_{s_1+1},...,a_{s_2})\in \mathcal{A}^{s_2-s_1+1}$ an admissible word where $s_1, s_2 \in \mathbb{Z} , s_1 < s_2$ and fix $s_1\le s\le s_2$. We write
	\begin{equation}\label{eq1}
	R(\alpha;s)=\bigcap_{m=s_1-s}^{s_2-s} \varphi^{-m}(R_{a_{m+s}}).
	\end{equation}

 Note that if $x\in R(\alpha;s)\cap \Lambda$ then the symbolic representation of $x$ is in the way $x=...a_{s_1}...a_{s-1};a_{s},a_{s+1}...a_{s_2}...$ where on the right of the ; is the $0$-th position. 
 
 We write $s^{(u)}(\alpha)$ for the unstable size  of $\alpha$, that is, the length  of the interval $I^u(\alpha)$
and the unstable scale of $\alpha$ is $r^{u}(\alpha)=\lfloor \log(1/s^{(u)}(\alpha))\rfloor$. Similarly, we write $s^{(s)}(\alpha)$ the stable size of $\alpha$ as being the length of $I^s(\alpha)$ and the stable scale of $\alpha$ is $r^{(s)}(\alpha) =\lfloor \log(1/s^{(s)}(\alpha))\rfloor.$

\begin{remark}\label{BDP}
In our context of $C^{1+\varepsilon}$-dynamically defined Cantor sets, we can relate the unstable and stable sizes of $\alpha$ to its length as a word in the alphabet $\mathcal{A}$ via the so-called \textbf{bounded distortion property} saying that there exists a constant $c_1 =c_1(\varphi,\Lambda) > 0$ such that
	\begin{equation}\label{bdp}
	e^{-c_1}\le \dfrac{|I^u(\alpha\beta)|}{|I^u(\alpha)||I^u(\beta)|}\le e^{c_1} \ \mbox{and} \ e^{-c_1}\le \dfrac{|I^s((\alpha\beta)^t)|}{|I^s(\alpha^t)||I^s(\beta^t)|}\le e^{c_1}.
	\end{equation}

\end{remark}
\begin{remark}\label{simmetry}
In the context of horseshoes of $C^2$-conservative diffeomorphisms, there is a constant $c_2=c_2(\varphi,\Lambda)>0$  such that the stable and unstable sizes of any word $\alpha=(a_1,...,a_n)$ in the alphabet $\mathcal{A}$ satisfy
\begin{equation}\label{simmetry}
e^{-c_2}\le \dfrac{|I^u(\alpha)|}{|I^s(\alpha^t)|}\le e^{c_2}.
\end{equation}
Indeed, this happens because $\varphi^n$ maps the unstable rectangle
	$$R^u(\alpha)=\{x\in R_{a_0}: \varphi^i(x)\in R_{a_i}, 1\le i\le n \}$$ 
diffeomorphically onto the stable rectangle
		$$R^s(\alpha^t)=\{x\in R_{a_0}: \varphi^j(x)\in R_{a_{n-j}}, 1\le j\le n\},$$
$\varphi$ preserves areas, and	the areas of $R^u(\alpha)$ and $R^s(\alpha^t)$ are comparable to $|I^u(\alpha)|$ and $|I^s(\alpha^t)|$ up to multiplicative factors.
\end{remark}
We write $\alpha^{\ast}=(a_1,a_2,...,a_{n-1})$ if $\alpha=(a_1,a_2,...,a_n)$ and for $r\in \mathbb{N}$ we define the sets
$$P^{(u)}_r=\{\alpha\in \mathcal{A}^n \ \mbox{admissible}:  r^{(u)}(\alpha)\geq r \ \mbox{and} \ r^{(u)}(\alpha^{\ast})<r\}$$
similarly,
$$P^{(s)}_r=\{\alpha\in \mathcal{A}^n \ \mbox{admissible}:  r^{(s)}(\alpha)\ge r \ \mbox{and} \ r^{(s)}(\alpha^{\ast})<r\}.$$
We also define
$$\mathcal{C}_u(t,r)=\{\alpha\in P^{(u)}_r:  I^{u}(\alpha)\cap K^u_t\neq \emptyset\}$$
and
$$ \mathcal{C}_s(t,r)=\{\alpha\in P^{(s)}_r:  I^{s}(\alpha^t)\cap K^s_t\neq \emptyset\}$$
whose cardinalities are denoted $N_u(t,r):=|\mathcal{C}_u(t,r)|$ and $N_s(t,r):=|\mathcal{C}_s(t,r)|$.

In the article \cite{CMM16} the authors proved the following lemma:

\begin{lemma}
For each $t\in \mathbb{R}$, the sequences $\{N_u(t,r)\}_{r\in \mathbb{N}}$, and $\{N_s(t,r)\}_{r\in \mathbb{N}}$ are essentially submultiplicative, in the sense that, there exists a constant $c= c(\varphi, \Lambda)\in \mathbb{N}$ such that
	$$N_u(t,m+n)\le |\mathcal{A}|^{c}\cdot N_u(t,m)\cdot N_u(t,n)$$
and
	$$N_s(t,m+n)\le |\mathcal{A}|^{c}\cdot N_s(t,m)\cdot N_s(t,n).$$
\end{lemma}
From this lemma, we get that for each $t\in \mathbb{R}$ there exist the limits 
\begin{equation}\label{D_u}
    D_u(t)=\lim_{r\to \infty}\dfrac{\log N_u(t,r)}{r}=\inf\limits_{r\in\mathbb{N}}\dfrac{\log (|\mathcal{A}|^{c}N_u(t,r))}{r}
\end{equation}
	and	
 \begin{equation}\label{D_s}
     D_s(t)=\lim_{r\to \infty}\dfrac{\log N_s(t,r)}{r}=\inf\limits_{r\in\mathbb{N}}\dfrac{\log (|\mathcal{A}|^{c}N_s(t,r))}{r}.
 \end{equation}
	
It is easy to see that the numbers $D_u(t)$ and $D_s(t)$ are the box-counting dimension of $K^u_t$ and $K^s_t$ respectively. By proposition 2.6 in \cite{CMM16} we have that $t\mapsto D_u(t)$ and $t\mapsto D_s(t)$ are upper semicontinuous functions without any condition on $HD(\Lambda)$. We proceed to give a proof that $t\mapsto D_u(t)$ and $t\mapsto D_s(t)$ are also lower semicontinuous functions, however in order to do that we need to work with the correct set of functions: 

Let $r\geq2$. We define
$$\mathcal{R}_{\varphi,\Lambda}:=\{f\in C^r(S,\mathbb{R}): \nabla f(z)\neq 0,  \ \forall \ z\in \Lambda\}.$$ 
In others words, $\mathcal{R}_{\varphi,\Lambda}$ is the class of functions $C^r$, $f :S \rightarrow \mathbb{R}$ such that for every $z\in \Lambda$ either $\ D f(z)e^s_z\neq 0$  or$\ D f(z)e^u_z\neq 0$ where $e^s_z$ and $e^u_z$ are unit vectors in the stable and unstable directions of $T_zS$. We end this section with the following lemma:
\begin{lemma}
Given $r\ge 2$, $\mathcal{R}_{\varphi,\Lambda}$ is an open and dense subset of $C^r(S,\mathbb{R}).$
\end{lemma}
\begin{proof}

Consider the class $\mathcal{M}$ of the Morse functions, we know by the compacity of $S$ that $\mathcal{M}$ is dense in $C^r(S,\mathbb{R})$ and as a corollary of Morse's lemma  that every element of $\mathcal{M}$ has only finitely many critical points. Then since we have $\inte(\Lambda)=\emptyset$, given $g\in \mathcal{M}$ we can find $f\in \mathcal{R}_{\varphi,\Lambda}$, $C^r$-arbitrarily close to $g$ and this implies that $\mathcal{R}_{\varphi,\Lambda}$ is also $C^r$-dense. As $\mathcal{R}_{\varphi, \Lambda}$ is clearly open we have the result.

\end{proof}

\section{Critical windows and combinatorial lemmas}

To prove the theorems we need the following proposition, whose proof depends on the notion of \textbf{critical window} and some combinatorial lemmas related with.

\begin{proposition}\label{prop1}
Let $\varphi:S\rightarrow S$ be a $C^2$ diffeomorphism with a transitive horseshoe $\Lambda$. Fix $f\in \mathcal{R}_{\varphi,\Lambda}$ and $t\in \mathbb{R}$ such that $D_u(t)$, resp. $D_s(t)>0$. Then, for every $0<\eta<1$ there exists $\delta>0$ and a complete subshift 
$\Sigma(\mathcal{B}_u) \subset \Sigma \subset \mathcal{A}^{\mathbb{Z}}\  \mbox{, resp.}\  \Sigma(\mathcal{B}_s)\subset \Sigma\subset \mathcal{A}^{\mathbb{Z}}$, associated to a finite set $\mathcal{B}_u=\{\beta^{(u)}_1,\beta^{(u)}_2,...,\beta^{(u)}_m\}$, resp. $\mathcal{B}_s=\{\beta^{(s)}_1,\beta^{(s)}_2,...,\beta^{(s)}_n\}$, of finite sequences, such that
	$$\Sigma(\mathcal{B}_u) \subset \Sigma_{t-\delta} \ \mbox{, resp.} \ \Sigma(\mathcal{B}_s) \subset \Sigma_{t-\delta}$$
and 	
	$$HD(K^u(\Sigma(\mathcal{B}_u))>(1-\eta)D_u(t) \ \mbox{, resp.} \ HD(K^s(\Sigma(\mathcal{B}_s^t))>(1-\eta)D_s(t).$$	
	where $K^u(\Sigma(\mathcal{B}_u))$ and $K^s(\Sigma(\mathcal{B}_s^t))$ are the subsets of $K^u$ and $K^s$, consisting of points whose trajectory under $g_u$ and $g_s$, follows an itinerary obtained from the concatenation of words in the alphabets $\mathcal{B}_u$ and $\mathcal{B}_s^t$ respectively, where $\mathcal{B}_s^t$ is the alphabet whose words are the transposes of the words of the alphabet $\mathcal{B}_s$.   
\end{proposition}

\begin{remark}\label{remark0}
By symmetry it suffices to exhibit $\mathcal{B}_u$ satisfying the conclusion of Proposition \ref{prop1}.
\end{remark}

\begin{remark}\label{remark1}
If the horseshoe $\Lambda$, and therefore the shift $\Sigma$ are transitive but not mixing, as in remark \ref{remarkA}, then any complete subshift (as $\Sigma(\mathcal{B}_u)$ and $\Sigma(\mathcal{B}_s)$) should be contained in a single mixing component $\Lambda^j$ of $\Lambda$, and thus $K^u(\Sigma(\mathcal{B}_u))$ and $K^s(\Sigma(\mathcal{B}_s^t))$ are dynamically defined Cantor sets.
\end{remark}

Fix $f\in \mathcal{R}_{\varphi,\Lambda}$ and let $m_1$ big enough such that if $\alpha= (a_{-s_1},...,a_0,...,a_{s_2})$ is admissible with $s_1 , s_2 > m_1$ then either $R(\alpha;0) \cap f^{-1}(t)=\emptyset$ or $R(\alpha;0) \cap f^{-1}(t)$ is the graph of a differentiable map $f_s$ defined in a (closed) sub interval of $I^s(a_0,a_{-1},...,a_{-s_1})$ with values in $I^u(a_0,a_1...,a_{s_2})$ (case 1) or $R(\alpha;0) \cap f^{-1}(t)$ is the graph of a differentiable map $f_u$ defined in a sub interval of $I^u(a_0,a_{1},...,a_{s_2})$ with values in $I^s(a_0,a_{-1}...,a_{-s_1})$ (case 2). Note that here we used the implicit function theorem, that we are working with the coordinates of the stable and unstable foliation, and also that we can suppose that with the choice made for $m_1$ there exists a $\tilde{\delta} >0$ such that in case 1: $\abs{Df(z)e^u _z} > \tilde{\delta}$, $\forall z \in R(\alpha;0)$ and in case 2: $\abs{Df(z)e^s _z} > \tilde{\delta}$,  $\forall z \in R(\alpha;0)$.

Suppose that we are in case 1. By the mean value theorem and because $f(\cdot,f_s( \cdot ))=t$, we have for some constant $M$ that depends only on $f$
\begin{equation}\label{Im}
    \abs{Im(f_s)}\leq \frac{M}{\tilde{\delta}} \abs{I^s(a_0,a_{-1}...,a_{-s_1})} \leq M_1 ({\lambda_{2,s}})^{s_1}, 
\end{equation}
where $\lambda_{2,s}$ is the greatest modulus of eigenvalues of $Df$ in $\Lambda$ at the stable direction  and $M_1$ is a constant. Suppose then $s_1 \geq R_1$ where $R_1$ is big enough such that there is at most only one $q\in \mathcal{A}$ with $I^u(a_0,...,a_{m _1},q) \cap Im(f_s) \neq \emptyset$.
\begin{figure}[h]
\centering
\includegraphics[width=0.65 \textwidth]{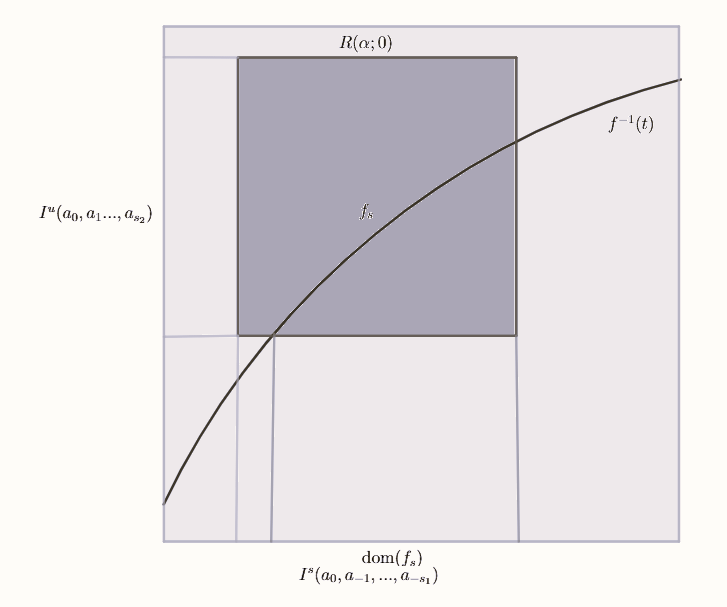}
\caption{Letters on the left of $\alpha$ determine part of the letters on the right.}
\end{figure}

Now assume that there exist $x_1, x_2 \in R(\alpha;0) \cap \Lambda$ such that $f(x_1)\leq t \leq f(x_2)$. Then we have two possibilities: 
If $a_{m_1 +1},...,{a}_{s_2} \in \mathcal{A}$ are unique with that property i.e. such that  there exist $x_1, x_2 \in R(\alpha;0) \cap \Lambda= R((a_{-s_1},...,a_{m_1},a_{m_1 +1},...,a_{s_2});0) \cap \Lambda$ with $f(x_1)\leq t \leq f(x_2)$, then by knowing $a_{-s_1},...,a_0,...,a_{m_1}$ we determine all the letters after $a_{m_1}$, i.e. the letters $a_{m_1 +1},...,a_{s_2}$.

If for some $j_0 \geq 1$ with $s_2 > j_0 + m_1$ there are $\tilde{a}_{m_1 + j_0+1},...,\tilde{a}_{s_2}, \tilde{b}_{m_1 + j_0+1},...,\tilde{b}_{s_2}\in \mathcal{A}$ with $\tilde{a}_{m_1 + j_0+1} \neq \tilde{b}_{m_1 + j_0+1} $ and $\tilde{x}_1, \tilde{x}_2 \in R((a_{-s_1},...,a_{m_1 + j_0 },\tilde{a}_{m_1 + j_0 + 1},...,\tilde{a}_{s_2});0) \cap \Lambda$, $\tilde{y}_1, \tilde{y}_2 \in R((a_{-s_1},...,a_{m_1 + j_0 },\tilde{b}_{m_1 + j_0 + 1},...,\tilde{b}_{s_2});0) \cap \Lambda$ with $f(\tilde{x}_1)\leq t \leq f(\tilde{x}_2)$ and $f(\tilde{y}_1)\leq t \leq f(\tilde{y}_2)$ then let $j_0$ minimal that satisfies that condition. So we have that depending on the relative positions of $I^u(a_0,...,a_{m_1 + j_0 },\tilde{a}_{m_1 + j_0+1})$, $I^u(a_0,...,a_{m_1 + j_0 },
\\ \tilde{b}_{m_1 + j_0+1} )$ and $Im(f_s)$ that $Im(f_s)$ contains an interval of the form $I^u(a_0,...,a_{m_1 + j_0 },q)$ with $q \in \mathcal{A}$ or contain a gap between two intervals of that form. In any case, by {\ref{Im}} and the mean value theorem again we have for some constant $C>0$ 
$$ C (\lambda_{2,u}^{-1})^{m_1 + j_0+1}  \leq M_1 ({\lambda_{2,s}})^{s_1},$$
where $\lambda_{2,u}$ is the greatest modulus of eigenvalues of $Df$ in $\Lambda$ at the unstable direction, and then for some $R>0$
\begin{equation}\label{Im2}
    (\lambda_{2,u}^{-1})^{j_0}  \leq R ({\lambda_{2,s}})^{s_1-m_1}.
\end{equation}

So by knowing the letters $a_{-m_1},...,a_0,...,a_{m_1}$ of $\alpha$, by \ref{Im2} the first $s_1-m_1$ letters determine 
$$ j_0 \geq \frac{\log(R)}{\log (\lambda_{2,u}^{-1})} + (s_1-m_1)\frac{\log(\lambda_{2,s})}{\log (\lambda_{2,u}^{-1})} >  (s_1 - m_1)\frac{\log(\lambda_{2,s})}{2 \log (\lambda_{2,u}^{-1})}$$
letters if $s_1 \geq R_2$ for some $R_2$.

Then either we determine $s_2-m_1$ letters or at least  
$$\left \lceil (s_1-m_1)\frac{\log(\lambda_{2,s})}{2 \log (\lambda_{2,u}^{-1})}\right \rceil \textrm \\\ \mbox{
letters if} \\\ s_1 \geq \max \{R_1, R_2\}.$$

If we are in case 2, we see analogously that are determine $s_1-m_1$ letters or at least
$$\left \lceil (s_2-m_1)\frac{\log(\lambda_{1,u}^{-1})}{2 \log (\lambda_{1,s})}\right \rceil \textrm \\\ \mbox{
letters if} \\\ s_2 \geq \max \{\tilde{R}_1, \tilde{R}_2\},$$
where $\lambda_{1,s}$, $\lambda_{1,u}$  are the smallest modulus of eigenvalues of $Df$ in $\Lambda$ at the stable and unstable direction respectivel and  $\tilde{R}_1, \tilde{R}_2$ are constants. 

Finally, in any case, if $s_1, s_2 \geq \max \{R_1,R_2,\tilde{R_1}, \tilde{R_2} \} := \tilde{R}$, $P$ letters at one side of the central block $(a_{-m_1},...,a_0,...,a_{m_1})$ either determine all the letters of the other side or at least 
\begin{equation}\label{m_0}
    \left \lceil \frac{1}{2}\min \left \{ \frac{\log(\lambda_{1,u}^{-1})}{ \log (\lambda_{1,s})}, \frac{\log(\lambda_{2,s})} { \log (\lambda_{2,u}^{-1})} \right \}P \right \rceil \geq \left \lceil \frac{1}{\frac{1}{\theta} + \tilde{R}}P   \right \rceil,
\end{equation}
where $\theta =\frac{1}{2} \min \left \{ \frac{\log(\lambda_{1,u}^{-1})}{ \log (\lambda_{1,s})}, \frac{\log(\lambda_{2,s})} { \log (\lambda_{2,u}^{-1})} \right\} < 1$.

Now, given $r\in \mathbb{N}$ define $\ell_1(r):=min \{|\beta|: \beta \in \mathcal{C}_u(t,r) \}$ and $\ell_2(r):=max \{|\beta|: \beta \in \mathcal{C}_u(t,r) \}$. For any word $\beta \in \mathcal{C}_u(t,r)$, as $r^{(u)}(\beta)\geq r$ and $r^{(u)}(\beta^*)< r$, we have for two constants $C_1,C_2>0$, with $\log C_1\notin \mathbb{Z}$ 
$$C_1 (\lambda_{2,u}^{-1})^{\abs{\beta}}  \leq  \abs{I^u(\beta)}\leq e ^{-r} < \abs{I^u(\beta^*)} \leq C_2 (\lambda_{1,u}^{-1})^{\abs{\beta}} $$
  then,
 $$\frac{-(r+\log C_1)}{\log(\lambda_{2,u}^{-1})}\leq \abs{\beta}<\frac{-(r+\log C_2)}{\log(\lambda_{1,u}^{-1})} $$
so, applying this to the words in $\mathcal{C}_u(t,r)$ that realize $\ell_1(r)$ and $\ell_2(r)$, we conclude that 
$$\frac{\ell_2(r)}{\ell_1(r)}\leq \frac{\log(\lambda_{2,u}^{-1})}{\log(\lambda_{1,u}^{-1})}\frac{(r+\log C_2)}{(r+\log C_1)}.$$
That is, $\{ \frac{\ell_2(r)}{\ell_1(r)} \}_{r\in \mathbb{N}}$ is bounded and then we can define 
  $$ m_0:=3\left \lceil \left( \frac{1}{\theta} + \tilde{R} \right) .\left( \sup \limits_{r\in \mathbb{N}}\frac{\ell_2(r)}{\ell_1(r)}\right)  \right \rceil.$$ 

In order to prove the proposition, let us begin by taking $\tau=\eta/(100(2m_0+3)^2)$
and $r_0=r_0(\varphi,f,\eta,t)\in \mathbb{N}$ large so that $\ell_1(r_0)\geq m_1$ and 
 \begin{equation}\label{eq r0}
\left | \dfrac{\log N_u(t,r)}{r} - D_u(t)\right | < \dfrac{\tau}{2}D_u(t), \ \forall \ r\ge r_0, \ r\in \mathbb{N}
\end{equation}  
also call $\mathcal{B}_0=\mathcal{C}_u(t,r_0)$, $N_0=N_u(t,r_0)$.

Consider $\beta=\beta_{k_1}\beta_{k_2}...\beta_{k_ \ell}=a_1...a_p \in \mathcal{A}^{p}, \ \beta_{k_i}\in \mathcal{B}_0, \ 1\le i\le \ell$. We say that $n\in \{1,...,p\}$ is the $n$-th position of $\beta$; if $\beta_{k_i}\in \mathcal{A}^{n_{k_{i}}} $ we write $|\beta_{k_i}|=n_{k_i}$ for its length and $P(\beta_{k_i})=\{1,2,...,n_{k_i}\}$ for its set of positions as a word in the alphabet $\mathcal{A}$ and given $s \in P(\beta_{k_i})$ we call  $P(\beta,k_i;s)=n_{k_1}+...+n_{k_{i-1}}+s$ the position in $\beta$ of the position $s$ of $\beta_{k_i}$.

For the next definition, set
$$\tilde{\mathcal{B}}=\tilde{\mathcal{B}_u}:=\{\beta=\beta_1\beta_2...\beta_k: \ \beta_j\in \mathcal{B}_0 \\\ \forall 1\leq j \leq k\ \mbox{and} \ I^u(\beta)\cap K^u_t\neq \emptyset\}$$
where $k=4(2m_0+3)N_0^{2m_0+3}\cdot \left \lceil 2/\tau \right \rceil.$

\begin{definition}\label{d2}
Let $\beta=\beta_1\beta_2...\beta_{k}, \ \beta_r\in \mathcal{B}_0, \ 1\leq r\leq k$ an element of $\tilde{\mathcal{B}}$. We say that $(i,j)$ is a critical window for $\beta$ if $j-i$ is even, $j-i\ge 2m_0+2$ and there is $n\in P(\beta_{(j+i)/2 })$ such that if $\tilde{\eta}=\beta_i...\beta_j=a_1...a_{|\tilde{\eta}|}$ there are 
$x_1,x_2\in R(\tilde{\eta};P(\eta,(j+i)/2;n))\cap \Lambda \ \ \mbox{with} \ \ f(x_1)\le t \le f(x_2).$ We call $r=(j-i)/2$ the radius of the critical window.

\end{definition}

\begin{remark}\label{cw}

By \ref{m_0} we have that in this situation we determine all the $r-1$ blocks of words of ${\mathcal{B}}_0$ that are on the left or on the right of $\beta_r, \beta_{r+1}, \beta_{r+2}$ or we determine 
$$j_0 \geq \left \lceil \frac{1}{\frac{1}{\theta} + \tilde{R}} \ell_1(r_0)(r-1)  \right \rceil  \geq \left \lceil \frac{1}{(\frac{1}{\theta} + \tilde{R})\frac{\ell_2(r_0)}{\ell_1(r_0)}} \ell_2(r_0)(r-1)  \right \rceil \geq  \left \lceil \frac{1}{m_0/3} \ell_2(r_0)(r-1)  \right \rceil$$
letters before the position $P(\tilde{\eta},r + 1;n)- m_1$ or after the position $P(\tilde{\eta},r + 1;n)+ m_1$ of $\tilde{\eta}$ and then 
$$ \left \lfloor \frac{r-1}{m_0/3} \right \rfloor-2 \geq  \left \lfloor \frac{r-1}{m_0} \right \rfloor$$
blocks at one side of $\beta_r, \beta_{r+1}, \beta_{r+2}$ are determined in any case.
\end{remark}
Given the pair $(i,j)$ we write $\overline{[i,j]}$ for the set $\{i,i+1,...,j\}.$ Moreover, if $\beta=\beta_1\beta_2...\beta_{k}, \ \beta_r\in \mathcal{B}_0, \ 1\le r\le k$ we put 
$$C(\beta)=\{1\le s\le k\ : \ \exists \ (i,j) \ \mbox{critical window of} \ \beta \ \mbox{and} \ s\in \overline{[i,j]}\}.$$
In other words, $C(\beta)$ is the set of positions that are ``contained" in a critical window.

Now we want to estimate the cardinality of the set
	$$\mathcal{E}=\left \{\beta=\beta_1...\beta_k \in \tilde{\mathcal{B}}: |C(\beta)|<\dfrac{k}{5(2m_0+3)}\right \}.$$
But first, we do that for the set $\tilde{\mathcal{B}}$.
\begin{lemma} We have
 $|\tilde{\mathcal{B}}|>2N_0^{(1-\tau)k}$.	
\end{lemma}
\begin{proof}
Given $\beta=\beta_1\beta_2...\beta_k \in \tilde{\mathcal{B}} $ we have by the inequality \ref{bdp} that 
	$$|I^u(\beta)|\le \prod_{i=1}^k e^{c_1}|I^u(\beta_i)|.$$

By definition, for every $i=1,2,...,k$, as $\beta_i\in \mathcal{C}_u(t,r_0)$
	$$r^u(\beta_i)=\left \lfloor \log \dfrac{1}{|I^u(\beta_i)|}\right \rfloor \ge r_0$$	
and then	
	$$|I^u(\beta)|\le \prod_{i=1}^k e^{c_1}|I^u(\beta_i)|\le e^{-k(r_0-c_1)}.$$
This implies that $\{I^u(\beta):\beta\in \tilde{\mathcal{B}}\}$ is a cover of $K^u_t$ by intervals of unstable-size $\ge k(r_0-c_1).$ In particular, writing $\beta=(b_1b_2...b_{n(k)})$ we have a surjective projection $(b_1b_2...b_{n(k)})\mapsto (b_1b_2...b_j)\in \mathcal{C}_u(t,k(r_0-c_1))$ where 
	$$j=\min\{1\le i\le n(k): r^u(b_1b_2...b_i))\ge k(r_0-c_1)\}.$$
We can take $r_0$ large enough such that $k(r_0-c_1)>r_0$ and then, by	(\ref{D_u}) 
	$$|\mathcal{C}_u(t,k(r_0-c_1))|=N_u(t,k(r_0-c_1))>\dfrac{1}{|\mathcal{A}|^{c}}e^{(k(r_0-c_1)D_u(t))}.$$
In particular,
 	$$|\tilde {\mathcal{B}}|\ge \dfrac{1}{|\mathcal{A}|^{c}}e^{(k(r_0-c_1)D_u(t))}>2e^{k(r_0-2c_1)D_u(t)},$$
because $k$ is large for $r_0$ large and $D_u(t)>0$. Then
	$$|\tilde {\mathcal{B}}|>2e^{(1-\tau/2)kr_0D_u(t)}>2e^{(1-\tau)(1+\tau/2)kr_0D_u(t)}>2 N_0^{(1-\tau)k},$$
because $N_0<e^{(1+\tau/2)r_0D_u(t)}$ by (\ref{eq r0}).	 		
\end{proof}
Using the above lemma we have the following:
\begin{lemma}\label{lem1}
One has $|\mathcal{E}|>N_0^{(1-\tau)k}.$
\end{lemma}
\begin{proof}
Remember the elementary fact that, given a finite family of intervals, there is a subfamily of disjoint intervals whose sum of lengths is at least half of the measure of the union of the intervals of the original family. If for $\beta \in \tilde {\mathcal B}$, $|C(\beta)| \geq \frac{k}{5(2m_0 +3)}$ then applying the above fact for the family of intervals $[i,j+1)$ with $(i,j)$ a critical window of $\beta$ there exist a family $\{ (i_x,j_x)\}_{x \in \mathcal{X}}$ of critical windows of $\beta$ such that $\overline{[i_x,j_x]} \cap \overline{[i_y,j_y]}=\emptyset $ if $x,y \in \mathcal{X}$ with $x \neq y$ and $\frac{k}{10(2m_0+3)} \leq |\bigcup_{x \in  \mathcal{X}} \overline{[i_x,j_x]}|:=M_\mathcal{X}$. Set $r_x=(j_x-i_x)/2$ for $x\in \mathcal{X}$. We observe that if $|\mathcal{X}| \leq \frac{M_\mathcal{X}}{2(2m_0+3)}$, then
\begin{eqnarray*}
\sum_{x\in \mathcal{X}} \left \lfloor \frac{r_x-1}{m_0} \right \rfloor &=& \frac{M_\mathcal{X}}{2m_0}-\frac{3|\mathcal{X}|}{2m_0}+\sum_{x\in \mathcal{X}} \left( \left \lfloor \frac{r_x-1}{m_0} \right \rfloor- \frac{r_x-1}{m_0} \right) \\ &=& \frac{M_\mathcal{X}}{2m_0}-\frac{3|\mathcal{X}|}{2m_0}-\sum_{x\in \mathcal{X}} \left( \frac{r_x-1}{m_0} - \left \lfloor \frac{r_x-1}{m_0} \right \rfloor \right)  \\ &\geq& \frac{M_\mathcal{X}}{2m_0}- \left(1+\frac{3}{2m_0} \right)|\mathcal{X}| \quad( \textrm{since } 0\leq \frac{r_x-1}{m_0} - \left \lfloor \frac{r_x-1}{m_0} \right \rfloor <1)\\ &\geq& \frac{M_\mathcal{X}}{2m_0}-\frac{M_\mathcal{X}}{4m_0}=\frac{M_\mathcal{X}}{4m_0} \geq \frac{k}{40m_0(2m_0+3)} \geq \frac{k}{20(2m_0+3)^2}
\end{eqnarray*}
and if $|\mathcal{X}| > \frac{M_\mathcal{X}}{2(2m_0+3)}$, then
 
$$ \sum_{x\in \mathcal{X}} \left \lfloor \frac{r_x-1}{m_0} \right \rfloor \geq \sum_{x\in \mathcal{X}} 1 = |\mathcal{X}| > \frac{M_\mathcal{X}}{2(2m_0+3)} \geq \frac{k}{20(2m_0+3)^2}.$$

In any case 
\begin{eqnarray*}
\prod_{x\in \mathcal{X}} N_0^{2r_x+1 -\lfloor (r_x-1)/m_0\rfloor} \cdot N_0^{k-M_\mathcal{X}}  &=&
    N_0^{M_\mathcal{X} - \sum_{x \in \mathcal{X}}\lfloor (r_x-1)/m_0\rfloor} \cdot N_0^{k-M_\mathcal{X}}\\&=& N_0^{k - \sum_{x \in \mathcal{X}}\lfloor (r_x-1)/m_0\rfloor}  \leq N_0^{(1-1/(20(2m_0+3)^2)k}.
\end{eqnarray*}

Then, using this and remark \ref{cw} we have:
\begin{equation} \label{tutu}
    |\tilde{\mathcal{B}}\setminus \mathcal{E}| \leq 2^k \cdot 2^k \cdot N_0^{(1-1/20(2m_0+3)^2)k}.
\end{equation}

Since for our choices of $r_0$, $N_0$, $k$ large enough 
 and $\tau$ sufficiently small we have $2^{2k}\cdot N_0^{(1-1/20(2m_0+3)^2 )k}<N_0^{(1-\tau)k}$ it follows from \ref{tutu} that:
$$|\mathcal{E}|=|\tilde{\mathcal{B}}|-|\tilde{\mathcal{B}}\setminus \mathcal{E}|>2N_0^{(1-\tau)k}-2^{2k}\cdot N_0^{(1-1/20(2m_0+3)^2 )k}>N_0^{(1-\tau)k}.$$	
This completes the proof of the lemma.
\end{proof}
Our next lemma shows that among the words $\beta \in \mathcal{E}$ we have several words which share the same positions which do not belong to $C(\beta)$ and the same words of  $\mathcal{B}_0$ appearing in these positions.  
\begin{lemma}\label{lem2}
 There are $3N_0^{2m_0+3}$ words $(\tilde{\beta}_{s_i-m_0-1},...,\tilde{\beta}_{s_i},...,\tilde{\beta}_{s_i+m_0+1})\in \mathcal{B}^{2m_0+3}_0$, with $s_i\in \{ m_0+2,m_0+3,....,k-m_0-1\}$, and $1\le i\le 3N_0^{2m_0+3}$,  such that
 	$$s_{i+1}-s_i\ge (2m_0+3)\left \lceil{ \dfrac{2}{\tau}}\right \rceil \ \mbox{for} \  1\le i< 3N_0^{2m_0+3} $$
and the set
	\begin{align*}X=\{\beta=\beta_1\beta_2...\beta_k\in \mathcal{E}: (\beta_{s_i-m_0-1},...,\beta_{s_i},...,\beta_{s_i+m_0+1})=(\tilde{\beta}_{s_i-m_0-1},...,\tilde{\beta}_{s_i},...,\tilde{\beta}_{s_i+m_0+1}), \\ \{s_i-m_0-1,...,s_i,...,s_i+m_0+1\}\cap C(\beta)=\emptyset,
	 1 \le i\le 3N_0^{2m_0+3}\}
	 \end{align*}	
has cardinality bigger than $N_0^{(1-2\tau)k}.$	
\end{lemma}

\begin{proof}
Given $\beta=\beta_1\beta_2...\beta_k\in \mathcal{E}$ we can find $W=\left \lceil \frac{4k}{5(2m_0+3)} \right \rceil$ indices $i_1<i_2<...<i_{W}$ with $i_p\in \{m_0+2,m_0+3,...k-m_0-1\}, \ \forall p=1,2,...,W$ such that
\begin{itemize}

\item $i_{p+1}-i_p\ge (2m_0+3), \ p=1,2,...,W-1$

\item $\cup_{p=1}^{W}\{i_p-m_0-1,...,i_p,...,i_p +m_0+1\}\cap C(\beta)=\emptyset.$

\end{itemize}
 We remember that $k=4(2m_0+3)N_0^{2m_0+3}\cdot \left \lceil 2/\tau \right \rceil$ and since $3N_0^{2m_0+3}\lceil 2/\tau\rceil < (16/5)N_0^{2m_0+3}\\ \lceil 2/\tau\rceil \leq W$ we can write $j_m=i_{m\lceil 2/\tau\rceil}$ with $1\le m\le 3N_0^{2m_0+3}$. Then for $1\le m< 3N_0^{2m_0+3}$, $j_{m+1}-j_m\ge (2m_0+3)\lceil 2/\tau\rceil$ and for $1\le m \leq 3N_0^{2m_0+3}$ 
	$$\{j_m-m_0-1,...,j_m,...,j_m+m_0+1\}\cap C(\beta)=\emptyset \   .$$
 
	Note that
	\begin{itemize}
	\item The number of possibilities for $(j_1,...,j_{3N_0^{2m_0+3}})$ is smaller than $2^k$	
	\item For $(j_1,...,j_{3N_0^{2m_0+3}})$ fixed, the number of possibilities for  $(\beta_{j_{i}-m_0-1},...,\beta_{j_{i}},...,\\ \beta_{j_{i} +m_0+1})$ is at most $N^{3(2m_0+3)N^{2m_0+3}_0}_0.$
	
	\end{itemize}
	Then we can choose $m_0+1<s_1<s_2<...<s_{3N_0^{2m_0+3}}<k-m_0$ with $s_{i+1}-s_i\ge (2m_0+3)\left \lceil{ 2/\tau}\right \rceil$ and strings $(\tilde{\beta}_{s_i-m_0-1},...,\tilde{\beta}_{s_i},...,\tilde{\beta}_{s_i+m_0+1}) \in \mathcal{B}^{2m_0+3}_0, \ 1\le i\le 3N_0^{2m_0+3}$ such that the set
		\begin{align*}X=\{\beta=\beta_1\beta_2...\beta_k\in \mathcal{E}: (\beta_{s_i-m_0-1},...,\beta_{s_i},...,\beta_{s_i+m_0+1})=(\tilde{\beta}_{s_i-m_0-1},...,\tilde{\beta}_{s_i},...,\tilde{\beta}_{s_i+m_0+1}), \\ \{s_i-m_0-1,...,s_i,...,s_i+m_0+1\}\cap C(\beta)=\emptyset,
	 1 \le i\le 3N_0^{2m_0+3}\}
	 \end{align*}	
	 has cardinality
	 $$|X|\geq \dfrac{|\mathcal{E}|}{2^k\cdot N_0^{3(2m_0+3)N_0^{2m_0+3}}}.$$
	 But, $|\mathcal{E}|>N_0^{(1-\tau)k}$ and $2^k\cdot N_0^{3(2m_0+3)N_0^{2m_0+3}}<N^{\tau k}_0.$ Therefore,
	 $$|X|>\dfrac{|\mathcal{E}|}{2^k\cdot N_0^{3(2m_0+3)N_0^{2m_0+3}}}>N_0^{(1-2\tau)k}.$$

\end{proof}

Our third combinatorial lemma states that it is possible to cut words of
the subset $X$ provided by Lemma \ref{lem2} at certain positions in such a way that
one obtains a set $\mathcal{B}_u$ with non-neglectible cardinality.

For every $1\leq p< q \leq 3N_0^{2m_0+3}$ we denote $\pi_{p,q}:X\rightarrow \mathcal{B}_0^{s_q-s_p}$ the projection $\pi_{p,q}(\beta)=(\beta_{s_p+1},\beta_{s_p+2},...,\beta_{s_q}),$ if $\beta=\beta_1\beta_2...\beta_k$.
\begin{lemma}\label{lem3}
There are $1\leq p_0<q_0\leq 3N_0^{2m_0+3}$ such that 
\begin{enumerate}
\item[i)] $(\tilde{\beta}_{s_{p_0}-m_0-1},...,\tilde{\beta}_{s_{p_0}},...,\tilde{\beta}_{s_{p_0}+m_0+1})=(\tilde{\beta}_{s_{q_0}-m_0-1},...,\tilde{\beta}_{s_{q_0}},...,\tilde{\beta}_{s_{q_0}+m_0+1})$
\item[ii)] $ |\pi_{p_0,q_0}(X)|>N_0^{(1-10\tau)(s_{q_0}-s_{p_0})}$

\end{enumerate}
\end{lemma}

\begin{proof}
Consider $\mathcal{T}$ the set of pairs $(p,q)$ such that $1\le p<q\le 3N_0^{2m_0+3}$ and $\abs{\pi_{p,q}(X)}\le N_0^{(1-10\tau)(s_q-s_p)}$. For each pair in $\mathcal{T}$ we exclude from the set $\overline{[1,3N_0^{2m_0+3}]}$ the indices $j\in \overline{[p,q-1]}$.

\textbf{Claim:} The set $\mathcal{Z}=\bigcup_{(p,q)\in \mathcal{T}} \overline{[p,q-1]}$ has cardinality smaller than $2N_0^{2m_0+3}.$ Using the same observation given in Lemma \ref{lem1} we can find a subset $\tilde{\mathcal{T}}$ of $\mathcal{T}$ such that  $\overline{[p,q-1]}\cap \overline{[\tilde{p},\tilde{q}-1]}=\emptyset,$ for every $(p,q),(\tilde{p},\tilde{q})\in \tilde{\mathcal{T}}$ with $(p,q)\neq (\tilde{p},\tilde{q})$ and
	$$\sum_{(p,q)\in \tilde{\mathcal{T}}}(q-p)\ge \dfrac{1}{2}|\mathcal{Z}|.$$
	
 Suppose that $|\mathcal{Z}|\ge 2N_0^{2m_0+3}$. Since the sequence $s_1<s_2<...<s_{3N_0^{2m_0+3}}$ given by Lemma \ref{lem2} is such that $s_{i+1}-s_i\geq (2m_0+3)\lceil 2/\tau \rceil$ we have
	\begin{equation}\label{eq6}
	\sum_{(p,q)\in \tilde{\mathcal{T}}}(s_{q}-s_{p})\ge  (2m_0+3)\lceil 2/\tau \rceil \sum_{(p,q)\in \tilde{\mathcal{T}}} (q-p) \ge (2m_0+3)\lceil 2/ \tau \rceil N_0^{2m_0+3}.
	\end{equation}
	On the other hand, since $\abs{\pi_{p,q}(X)}\le N_0^{(1-10\tau)(s_q-s_p)}$ we have
	$$|X|<N_0^{(1-10\tau)\sum_{(p,q)\in \tilde{\mathcal{T}}}(s_q-s_p)}\cdot N_0^{k-\sum_{(p,q)\in \tilde{\mathcal{T}}}(s_q-s_p)}$$
	
 Using \ref{eq6} we have
	\begin{equation}\label{eq7}
	|X|< N_0^{(1-10)\tau((2m_0+3)\lceil 2/ \tau \rceil N_0^{2m_0+3})} \cdot N_0^{k-(2m_0+3)\lceil 2/ \tau \rceil N_0^{2m_0+3})}              =N_0^{k-10\tau((2m_0+3)\lceil 2/ \tau \rceil N_0^{2m_0+3})}.
	\end{equation}
	By Lemma \ref{lem2} we know that $|X|>N_0^{(1-2\tau)k}$. Using that and the inequality (\ref{eq7}) we must have
	$$(1-2\tau)k<k-10\tau((2m_0+3)\lceil 2/ \tau \rceil N_0^{2m_0+3})$$
that is,
	$$10(2m_0+3)\lceil 2/ \tau \rceil N_0^{2m_0+3}<2k.$$
Since, $2k=8(2m_0+3)\lceil 2/ \tau \rceil N_0^{2m_0+3}$	 we have a contradiction. This implies that $|\mathcal{Z}|<2N_0^{2m_0+3}$ which proves our claim.

Therefore, we do not exclude at least $N_0^{2m_0+3}+1$ indices. Since for each of that indices we have at most $N_0^{2m_0+3}$ possibilities for choose $(\tilde{\beta}_{s_i-m_0-1},...,\tilde{\beta}_{s_i},...,\tilde{\beta}_{s_i+m_0+1})$ (see lemma \ref{lem2}) we conclude that there are two indices $(p_0,q_0) \notin \mathcal{T}$ such that
	$$(\tilde{\beta}_{s_{p_0}-m_0-1},...,\tilde{\beta}_{s_{p_0}},...,\tilde{\beta}_{s_{p_0}+m_0+1})=(\tilde{\beta}_{s_{q_0}-m_0-1},...,\tilde{\beta}_{s_{q_0}},...,\tilde{\beta}_{s_{q_0}+m_0+1}).$$
	By definition of non-excluded index $ |\pi_{p_0,q_0}(X)|>N_0^{(1-10\tau)(s_{q_0}-s_{p_0})}$ as we wanted to see.
\end{proof}
    Take $\mathcal{B}_u:=\pi_{p_0,q_0}(X)$ were $p_0, q_0$ are given by the previous lemma.
	Note that $K^u(\Sigma(\mathcal{B}_u))$ is a $C^{1+\varepsilon}$-dynamically defined Cantor set associated to certain iterates of $g_u$ on the intervals $I^u(\alpha)$ with $\alpha \in \mathcal{B}_u$. In this case its Hausdorff dimension coincides with its box-counting dimension and then
	\begin{equation}\label{eq 1, lem3}HD(K^u(\Sigma(\mathcal{B}_u)))\ge \left (1-\dfrac{\tau}{2}\right) \dfrac{| \mathcal{B}_u|}{- \log (\min \limits_{\alpha\in      \mathcal{B}_u} |I^u(\alpha)|)}
	\end{equation}
for $r_0$ sufficiently large. By the item ii) of lemma \ref{lem3},  $|\mathcal{B}_u|>N_0^{(1-10\tau)(s_{q_0}-s_{p_0})}$. On the other hand, by the bounded distortion property (see \ref{bdp}), we have $|I^u(\alpha)|\ge e^{-(c_1+r_0)(s_{q_0}-s_{p_0})}$, for each $\alpha\in \mathcal{B}_u$. Using this and the inequality \ref{eq 1, lem3} we obtain
	\begin{equation}\label{eq 2, lem3}
	HD(K^u(\Sigma(\mathcal{B}_u)))\ge \dfrac{(1-\tau/2)(1-10\tau) \log N_0}{c_1+r_0}.
	\end{equation}

Since $N_0=N_u(t,r_0)$ satisfies (see \ref{eq r0})
	$$\left | \dfrac{\log N_0}{r_0} - D_u(t)\right | < \dfrac{\tau}{2}D_u(t)$$
	we have
	$$\log N_0>(1-\tau/2)r_0D_u(t).$$
Plugging this in the inequality \ref{eq 2, lem3} and using that $\tau=\eta/(100(2m_0+3)^2)$ we have
	$$HD(K^u(\Sigma(\mathcal{B}_u)))>\dfrac{(1-10\tau)(1-\tau/2)^2r_0}{r_0+c_1}D_u(t)>(1-12\tau)D_u(t)>(1-\eta)D_u(t)$$
for $r_0=r_0(\eta)$ sufficiently large.	

At this point, we are ready to end the proof of Proposition \ref{prop1}.  

\begin{proof}
 We write by simplicity
$$y_1=\tilde{\beta}_{s_{p_0}+1}\tilde{\beta}_{s_{p_0}+2}...\tilde{\beta}_{s_{p_0}+m_0+1}=\tilde{\beta}_{s_{q_0}+1} \tilde{\beta}_{s_{q_0}+2} ...\tilde{\beta}_{s_{q_0}+m_0+1}$$
 and
$$y_2=\tilde{\beta}_{s_{p_0}-m_0-1}...\tilde{\beta}_{s_{p_0}-1}\tilde{\beta}_{s_{p_0}}=\tilde{\beta}_{s_{q_0}-m_0-1}...\tilde{\beta}_{s_{q_0}-1}\tilde{\beta}_{s_{q_0}}.$$

 It follows that any element in $\mathcal{B}_u$ has the form $y_1\beta_{s_{p_0}+m_0+2}...\beta_{s_{q_0}-m_0-2}y_2$, where $\beta_i\in \mathcal{B}_0$ for any $i=s_{p_0}+m_0+2,...,s_{q_0}-m_0-2$ and 
	 $$I^u(y_2y_1\beta_{s_{p_0}+m_0+2}...\beta_{s_{q_0}-m_0-2}y_2y_1)\cap K^u_t\neq \emptyset.$$
And then the elements of $\Sigma(\mathcal{B}_u)$ have the form
	$$x= \theta^{(1)}y_2;y_1\beta_{s_{p_0}+m_0+2}...\beta_{s_{q_0}-m_0-2}y_2y_1\theta^{(2)},$$ 
where ; indicates that the $0$-th position is the first position in $y_1$ and $y_1\theta^{(2)}\in \mathcal{A}^{\mathbb{N}}$ and $\theta^{(1)}y_2 \in \mathcal{A}^{\mathbb{Z}^-}.$	

We will see that $f(\sigma^\ell (x))<t$ ,  $ \forall \ell\in \mathbb{Z}$. Repair that in order to prove that, it is sufficient consider $0 \leq  \ell \leq \tilde{m}-1$ if 
     $$\alpha=y_1\beta_{s_{p_0}+m_0+2}...\beta_{s_{q_0}-m_0-2}y_2=a_1 \dots a_{\tilde{m}}.$$
Take then $j\in P(\beta_r)$ for some $r\in \{s_{p_0}+1,...,s_{q_0} \}$ and suppose that $f(\sigma^{P(\alpha,r; j)}(x))\geq t$. If $s_{q_0}-r \geq r-s_{p_0}-1$ let $\tilde{\eta}= \beta_{s_{p_0}-m_0}...\beta_r...\beta_{2r-s_{p_0}+m_0}$, then $x_1=\sigma^{P(\alpha,r; j)}(x)\in R(\tilde{\eta}; P(\tilde{\eta},r; j)) \cap \Lambda$ and since $I^u(y_2y_1\beta_{s_{p_0}+m_0+2}...\beta_{s_{q_0}-m_0-2}y_2y_1)\cap K^u_t\neq \emptyset$, by definition there are $\theta^{(3)}\in \mathcal{A}^{\mathbb{Z}_{-}}$ and $\theta^{(4)}\in \mathcal{A}^{\mathbb{N}}$ such that
	$$\theta^{(3)};y_2y_1\beta_{s_{p_0}+m_0+2}...\beta_{s_{q_0}-m_0-1}y_2y_1\theta^{(4)}\in \Sigma_t,$$
and then, there exists $x_2\in R(\tilde{\eta}; P(\tilde{\eta},r; j)) \cap \Lambda$ such that $f(x_2)\le t$. But this is a contradiction because remembering that $\mathcal{B}_u=\pi_{p_0,q_0}(X)$ then there is some $\beta \in X$ such that $(s_{p_0}-m_0,2r-s_{p_0}+m_0)$ is a critical window of $\beta$, because $2r-s_{q_0}+m_0-(s_{p_0} -m_0)=2r-2s_{p_0}+2m_0 \geq 2m_0+2$, and $s_{p_0} \in \overline{[s_{p_0}-m_0,2r-s_{p_0}+m_0]}$.

If $s_{q_0}-r < r-s_{p_0}-1$ the argument is similar. Therefore,  $f(\sigma^\ell (x))<t$ ,  $ \forall \ell\in \mathbb{Z}$ and since $\Sigma(\mathcal{B}_u) \subset \bigcup_{i=0}^{k\ell_2(r_0)} \sigma ^i(\Sigma(\mathcal{B}_u))$ where $\bigcup_{i=0}^{k\ell_2(r_0)} \sigma ^i(\Sigma(\mathcal{B}_u))=\bigcup_{i\in \mathbb{Z}} \sigma ^i(\Sigma(\mathcal{B}_u))$ is the compact set, formed by the orbits by $\sigma$ of elements of $\Sigma(\mathcal{B}_u)$, there exists $\delta>0$ such that
	$$\Sigma(\mathcal{B}_u)\subset \Sigma_{t-\delta}.$$	
	
\end{proof}

\begin{remark}
It is possible to show that if $\Sigma(\mathcal{B})\subset \Sigma\subset
\mathcal{A}^{\mathbb{Z}}$ is a complete subshift associated to a 
finite alphabet $\mathcal{B}$ of finite words on $\mathcal{A}$ then the set of the previous proof (as a subset of $\Lambda$) $\Lambda(\Sigma(\mathcal{B}))=\Pi^{-1}(\bigcup_{i\in \mathbb{Z}} \sigma ^i(\Sigma(\mathcal{B})))$ is a subhorseshoe of $\Lambda$.
\end{remark}

\section{Proofs of the Theorems 1 and 2}

We begin proving Theorem 1.
\begin{proof}
 The Proposition \ref{prop1} implies that
 $$D_u(t)\ge HD(K^u_t)\ge HD(K^u_{t-\delta})\ge HD(K^u(\Sigma(\mathcal{B}_u)))>(1-\eta)D_u(t).$$
 Since $\eta>0$ is arbitrary we have $D_u(t)=HD(K^u_t)=d_u(t).$ Moreover,
	$$(1-\eta)D_u(t)\le HD(K^u(\Sigma(\mathcal{B}_u)))\le HD(K^u_{t-\delta})=D_u(t-\delta)$$
that is, $t\mapsto D_u(t)$ is a lower semicontinuous function. Since by Proposition 2.6 in \cite{CMM16}, $t\mapsto D_u(t)$ is also upper semicontinuous, we have that $t\mapsto D_u(t)=d_u(t)$
is continuous.

Similarly, we have the equality $D_s(t)=d_s(t)$ and that $a\mapsto D_s(t)=d_s(t)$ is continuous, so we have proved theorem 1.	
\end{proof}

In the sequel, we will use the following result that follows from the spectral decomposition theorem and from \cite{M50}

\begin{proposition}\label{theorem M}
In the context of Theorem 2. There exists a residual subset $\tilde{\mathcal{U}}\subset \mathcal{U}$ with the property that for every subhorseshoe $\tilde{\Lambda}\subset \Lambda$ and any $f\in C^1(S,\mathbb{R})$ such that there exists some point in $\tilde{\Lambda}$ with its gradient not parallel neither the stable direction nor the unstable direction, one has 
	$$HD(f(\tilde{\Lambda}))=\min\{1,HD(\tilde{\Lambda})\}.$$
\end{proposition}
to prove the next proposition

\begin{proposition} \label{lagrange}
If $\tilde{\mathcal{U}}$ is as in the proposition \ref{theorem M} then for any $\varphi \in \tilde{\mathcal{U}}$, there exists a $C^r$-residual subset $\tilde{\mathcal{R}}_{\varphi,\Lambda}\subset \mathcal{R}_{\varphi,\Lambda}$ such that for every subhorseshoe $\widetilde{\Lambda}\subset\Lambda$ and any $f\in \tilde{\mathcal{R}}_{\varphi,\Lambda}$ one has 
$$\min\{1,HD(\tilde{\Lambda})\}=HD(\ell_{\varphi,f}(\widetilde{\Lambda}))=HD(m_{\varphi,f}(\widetilde{\Lambda})).$$
\end{proposition}

\begin{proof}
Following the ideas of the proof of the theorem 1 of \cite{MR2} we see that given a subhorseshoe $\widetilde{\Lambda}\subset\Lambda$, the set 
$$H_{\widetilde{\Lambda}}=\{f\in C^r(S,\mathbb{R}): \abs{M_{\widetilde{\Lambda}, f}}=1\ \text{and if}\ z\in M_{\widetilde{\Lambda}, f},\ Df_z(e^{s,u}_z)\neq 0  \}$$
is $C^r$- open and dense set, where $M_{\widetilde{\Lambda}, f}= \{x\in \widetilde{\Lambda}: \forall y\in \widetilde{\Lambda} \mbox{ , } f(x) \geq f(y) \}$. Take then  
$$\tilde{\mathcal{R}}_{\varphi,\Lambda}:= \bigcap \limits_{ \substack{\widetilde{\Lambda}\subset\Lambda\  \\ subhorseshoe}}H_{\widetilde{\Lambda}} \cap \mathcal{R}_{\varphi, \Lambda }.$$



In the mentioned paper is also proved that for any such subhorseshoe $\widetilde{\Lambda}\subset\Lambda$ and $f \in \tilde{\mathcal{R}}_{\varphi,\Lambda}$ if $x_M$ is the unique element where $f|_{\widetilde{\Lambda}}$ take its maximum value, then for any $\epsilon>0$ there exists some subhorseshoe $\widetilde{\Lambda}^{\epsilon}\subset \widetilde{\Lambda} \setminus \{ x_M\}$ with 
$$HD(\widetilde{\Lambda}^{\epsilon})\geq HD(\widetilde{\Lambda})(1-\epsilon)$$ 
and such that for some point $d\in \widetilde{\Lambda}^{\epsilon}$ there exists a local $C^{1}$-diffeomorphism $\tilde{A}$ defined in a neighborhood $U_{d}$ of $d$ such that 
$$f(\varphi^{j_0}(\tilde{A}(\tilde{\Lambda}_{j_0})))\subset \ell_{\varphi,f}(\widetilde{\Lambda}),$$
where $j_{0}$ is an integer and $\tilde{\Lambda}_{j_0}\subset \widetilde{\Lambda}^{\epsilon}$ has nonempty interior in $\widetilde{\Lambda}^{\epsilon}$ and then is such that $HD(\tilde{\Lambda}_{j_0})=HD(\widetilde{\Lambda}^{\epsilon})$. Moreover, it is proved also that $\dfrac{\partial \tilde{A}}{\partial e_{x}^{s,u}}\parallel e^{s,u}_{\tilde{A}(x)}$, for $x\in U_{d}\cap  \widetilde{\Lambda}^{\epsilon}$ and then, by construction, $\nabla (f\circ \varphi^{j_{0}} \circ \tilde{A})(x)\nparallel e_{x}^{s,u}$ for every $x \in \tilde{\Lambda}_{j_0}$.

Extending properly $f\circ \varphi^{j_{0}} \circ \tilde{A}$, and letting $\epsilon$ tends to $0$; it follows from this and proposition \ref{theorem M} that
$$ \min\{1,HD(\tilde{\Lambda})\}\leq HD(\ell_{\varphi,f}(\widetilde{\Lambda})).$$
And finally
$$\min\{1,HD(\tilde{\Lambda})\}\leq HD(\ell_{\varphi,f}(\widetilde{\Lambda})) \leq HD(m_{\varphi,f}(\widetilde{\Lambda})) \leq HD(f(\widetilde{\Lambda})) \leq \min\{1,HD(\tilde{\Lambda})\}.$$
As we wanted to see.
\end{proof}
Now we proceed with the proof of Theorem 2.
\begin{proof}
 First, note that as in \ref{eq 1, lem3} we have
	$$HD(K^s(\Sigma(\mathcal{B}^t_u)))>\left ( 1-\dfrac{\tau}{2} \right ) \dfrac{\log |\mathcal{B}_u|}{-\log (\min_{\alpha\in \mathcal{B}_u} |I^s(\alpha^t)|)},$$
where $\mathcal{B}^t_u$ is the alphabet whose words are the transposes of the words of the alphabet $\mathcal{B}_u$. Since $|I^s(\alpha^t)|$ is comparable to $|I^u(\alpha)|$, using the notation of the remark \ref{simmetry} and the calculations after \ref{eq 1, lem3} we have that for $r_0$ large
	$$D_s(t)\geq HD(K^s(\Sigma(\mathcal{B}^t_u))) \geq \dfrac{(1-10\tau)(1-\tau/2)^2r_0D_u(t)}{r_0+c_1+c_2}>(1-\eta)D_u(t).$$
Since $\eta>0$ is arbitrary we have $D_s(t)\geq D_u(t)$ and the other inequality is proved in a similar way. On the other hand, if we take $\varphi\in \tilde{\mathcal{U}}\subset \mathcal{U}$, $t\in \mathbb{R}$ such that $D_u(t)>0$ and $\eta>0$ we have
\begin{equation}\label{Lambda_t}
    2(1-\eta)D_u(t)=(1-\eta)(D_s(t)+D_u(t))\le HD(\Lambda(\Sigma(\mathcal{B}_u))),
\end{equation}	
where $\mathcal{B}_u$ comes from Proposition \ref{prop1}. By Proposition \ref{lagrange} it follows that
	$$\min\{1,HD(\Lambda(\Sigma(\mathcal{B}_u)))\}= HD(\ell_{\varphi,f}(\Lambda(\Sigma(\mathcal{B}_u))))$$
and then
\begin{eqnarray*}
\min\{1,2(1-\eta)D_u(t)\}&\leq& \min\{1,HD(\Lambda(\Sigma(\mathcal{B}_u)))\}= HD(\ell_{\varphi, f}(\Lambda(\Sigma(\mathcal{B}_u)))) \\ 
&\leq& HD(L_{\varphi,f}\cap(-\infty,t))\leq HD(M_{\varphi,f}\cap(-\infty, t)) \\ 
&\leq& HD(f(\Lambda_{t})) \leq \min\{1,HD(\Lambda_{t})\} \\  &\leq& \min\{1,2D_u(t)\}.
\end{eqnarray*}

Since $\eta>0$ is arbitrary 		
	$$\min\{1,2D_u(t)\}=L(t)=M(t).$$

Finally, using one more time \ref{Lambda_t}, we also obtain
\begin{eqnarray*}
2(1-\eta)D_u(t) &\leq& HD(\Lambda(\Sigma(\mathcal{B}_u))) 
 \leq HD(\Lambda_t) \leq 2D_u(t),
\end{eqnarray*}
because $\eta>0$ is arbitrary, this proves that $HD(\Lambda_t)=2D_u(t)$.
\end{proof}

\begin{remark}
    The equality $HD(\Lambda_t)=2D_u(t)$ in the last proof, in fact, doesn't need any generic condition on $\varphi$.
\end{remark}

Now we want to prove that the conclusions of proposition \ref{lagrange} hold not only for subhorseshoes, but also for sets of the form $\Lambda_t$ for $t\in \mathbb{R}$. In order to do that, we have the following lemma

\begin{lemma}\label{L}
For every $t\in\mathbb{R}$ we have 
$$L(t)=\sup \limits_{s <t} HD(\ell_{\varphi,f}(\Lambda_s))=\lim \limits_{s \to\ t^-} HD(\ell_{\varphi,f}(\Lambda_s))$$
and 
$$M(t)=\sup \limits_{s <t} HD(m_{\varphi,f}(\Lambda_s))=\lim \limits_{s \to\ t^-} HD(m_{\varphi,f}(\Lambda_s)).$$
\end{lemma}

\begin{proof}
  Let $x\in \Lambda$ with $\ell_{\varphi,f}(x)=\eta < t$, then there exist a sequence $\{ n_k \}_{k\in\mathbb{N}}$ such that $\lim \limits_{k \to \infinity} f(\varphi^{n_k}(x))=\eta$. By compactness, without loss of generality, we can also suppose that $\lim \limits_{k \to \infinity} \varphi^{n_k}(x)=y$ for some $y\in \Lambda$ and so that $f(y)=\eta$. 

We affirm that $m_{\varphi, f}(y)=\eta$: in other case we would have for some $\tilde{k} \in \mathbb{Z}$ and $r \in \mathbb{R}$, $f(\varphi^{\tilde{k}}(y))>r>\eta$ and then for $k$ big enough by continuity $f(\varphi^{\tilde{k}+n_k}(x))>\eta$ that contradicts the definition of $\eta$. 
Now, take $\tilde{N}$ big enough such that if for two elements $a, b\in \Lambda$ their kneading sequences coincide in the central block (centered at the zero position) of size $2\tilde{N}+1$ then $\abs{f(a)-f(b)}<(t-\eta)/4$ and $N$ big enough such that for $k\geq N$ one has $f(\varphi^{k}(x))<\eta+(t-\eta)/4$ and the kneading sequences of $\varphi^{n_k}(x)$ and $y$ coincide in the central block of size $2\tilde{N}+1$. Suppose $\Pi(x)= (x_n)_{n\in \mathbb{Z}}$ and $\Pi(y)= (y_n)_{n\in \mathbb{Z}}$, then the point
$$\tilde{y}=\Pi^{-1}(\dots,y_{-n},\dots,y_{-1},x_{n_N},x_{n_N+1},\dots)$$
satisfies $\tilde{y}\in \Lambda_{(t+\eta)/2}$ and $\ell_{\varphi,f}(\tilde{y})=\eta.$ Therefore, we conclude that 
$$\ell_{\varphi,f}( \ell_{\varphi,f}^{-1}(-\infinity,t))=\ell_{\varphi,f}( \{x\in \Lambda: \ell_{\varphi,f}(x)<t \}) \subset \bigcup \limits_{s<t} \ell_{\varphi,f}(\Lambda_s)$$ 
and as for $s<t$, \ $\Lambda_s \subset \ell_{\varphi,f}^{-1}(-\infinity,t)$, the other inclusion also holds and we have the equality 
$$\ell_{\varphi,f}( \ell_{\varphi,f}^{-1}(-\infinity,t))= \bigcup \limits_{s<t} \ell_{\varphi,f}(\Lambda_s).$$
From this follows the result
\begin{eqnarray*}
L(t)&=&HD(\ell_{\varphi,f}( \ell_{\varphi,f}^{-1}(-\infinity,t)))=HD(\bigcup \limits_{s<t} \ell_{\varphi,f}(\Lambda_s))= HD(\bigcup \limits_{n\in \mathbb{N}} \ell_{\varphi,f}(\Lambda_{t-1/n}))\\ &=&\sup \limits_{n\in \mathbb{N}} HD(\ell_{\varphi,f}(\Lambda_{t-1/n}))=\sup \limits_{s<t} HD(\ell_{\varphi,f}(\Lambda_s)).
\end{eqnarray*}

For the second identity, as before 
$$M(t)=HD(M_{\varphi,f}(\Lambda)\cap (-\infty,t))=\sup \limits_{s<t}HD(M_{\varphi,f}(\Lambda)\cap (-\infty,s])= \sup \limits_{s<t} HD(m_{\varphi,f}(\Lambda_s)).$$  
\end{proof}

\begin{corollary}
  For any $\varphi \in \tilde{\mathcal{U}}$, $f\in \tilde{\mathcal{R}}_{\varphi,\Lambda}$ and $t\in \mathbb{R}$
$$\min\{1,HD(\Lambda_t)\}=HD(\ell_{\varphi,f}(\Lambda_t))=HD(m_{\varphi,f}(\Lambda_t)).$$
  \end{corollary}

  \begin{proof}
      This is a direct consequence of lemma \ref{L} and  theorem \ref{main2}. Indeed, for $\delta>0$
      $$L(t-\delta)\leq HD(\ell_{\varphi,f}(\Lambda_t))\leq L(t+\delta)$$
      letting $\delta$ tends to $0$ we have $L(t)=HD(\ell_{\varphi,f}(\Lambda_t)).$ Analogously we have $M(t)=HD(m_{\varphi,f}(\Lambda_t))$ and from theorem \ref{main2}, $L(t)=M(t)=\min\{1,HD(\Lambda_{t})\}.$
  \end{proof}

We end this work by giving another property of the map $L=M$

\begin{corollary}
For any $\varphi \in \tilde{\mathcal{U}}$ and $f\in \tilde{\mathcal{R}}_{\varphi,\Lambda}$ the map $L=M$ is not a Holder continuous function.
\end{corollary}

\begin{proof}

 By proposition \ref{lagrange}, for $\varphi \in \tilde{\mathcal{U}}$ and $f\in \tilde{\mathcal{R}}_{\varphi,\Lambda}$ one has 
 $$HD(L_{\varphi, f})=HD(\ell_{\varphi,f}(\Lambda))=\min\{1,HD(\Lambda)\}>0.$$
 As $L$ is continuous and $Im(L)=L(\mathbb{R})=[0,HD(L_{\varphi, f})]$, we can consider the point $c_{\varphi,f}$  where the map $L$ begins to be positive.

Suppose that $L$ is Holder continuous with exponent $\alpha>0$. Then there is $\epsilon >0$ such that $0<L(c_{\varphi,f}+\epsilon)<\alpha$ and being $L$ an $\alpha$-Holder function, one has 
\begin{eqnarray*}
1=HD([0,L(c_{\varphi,f}+\epsilon)])&=&HD(L(L_{\varphi,f}\cap(-\infty,c_{\varphi,f}+\epsilon))) \\ &<&\frac{1}{\alpha}.HD(L_{\varphi,f}\cap(-\infty,c_{\varphi,f}+\epsilon))\\ &=&\frac{1}{\alpha}L(c_{\varphi,f}+\epsilon) < 1,
\end{eqnarray*} 
which is an absurd.
\end{proof}


\begin{thebibliography}{10}



\bibitem{CF-89}
T.~W. Cusick and M.~E. Flahive.
\newblock {\em Markoff and Lagrange Spectra}, volume~30.
\newblock Mathematical Surveys and Monographs, Amer. Math. Soc., 1989.

\bibitem{CMM16}
A.~Cerqueira, C. ~Matheus, C. G.~Moreira.
\newblock Continuity of Hausdorff dimension across generic dynamical Lagrange
  and Markov spectra.
\newblock {\em Journal of Modern Dynamics}, 12:151--174, 2018.  \url{https://www.aimsciences.org/article/doi/10.3934/jmd.2018006}. 
\bibitem{LMMR}
D.~Lima, C. Matheus, C. G. Moreira and S. Roma\~na 
\newblock   {\em Classical and Dynamical Markov and Lagrange spectra},
\newblock{\em World Scientific}, 2020.

\bibitem{LM}
D.~Lima and C. G. Moreira.
\newblock  Phase transtitions on the Markov and Lagrange dynamical spectra.
\newblock{\em Annales de L'Institute Henri Poincar\'e (C), Analyse non-lin\'eaire}, Volume 38, Issue 5,
2021, Pages 1429-1459, ISSN 0294-1449,
\url{https://doi.org/10.1016/j.anihpc.2020.11.007.}

\bibitem{BK}
B. P. ~Kitchens.
\newblock {\em Symbolic Dynamics: One-sided, Two-sided and Countable State Markov Shifts},
\newblock Universitext, Springer, 1997.

\bibitem{M50}
C. G. Moreira.
\newblock Geometric properties of images of cartesian products of regular Cantor
sets by differentiable real maps,
\newblock {\em Mathematische Zeitschrift}, (2023) 303:3, \url{https://doi.org/10.1007/s00209-022-03151-z}.

\bibitem{MMan}
H. ~McCluskey and A. Manning.
\newblock Hausdorff dimension for horseshoes.
\newblock {\em Ergodic Theory and Dynamical Systems}, 3:251--260, 1983, \url{ https://doi.org/10.1017/S0143385700001966}.

\bibitem{M79}
A.~Markoff.
\newblock Sur les formes quadratiques binaires ind\'efinies.
\newblock {\em Math.Ann.}, 15:381--406, 1879, \url{https://eudml.org/doc/156934}.

\bibitem{M3}
C.~G. Moreira.
\newblock Geometric properties of Markov and Lagrange spectra.
\newblock {\em Annals of Math.}, 188: 145--170, 2018, \url{https://doi.org/10.4007/annals.2018.188.1.3}

\bibitem{MR2}


C.~G. Moreira and S.~Roma\~na.
\newblock On the Lagrange and Markov dynamical spectra.
\emph{Ergodic Theory and Dynamical Systems}, \textbf{37(5)}. 1570-1591, \url{https://doi.org/10.1017/etds.2015.121}.

\bibitem{CF}
T.~Cusick and M.~Flahive, 
\newblock{The Markoff and Lagrange spectra}, 
\emph{Mathematical Surveys and Monographs}, \textbf{30}. American Mathematical Society, Providence, RI, 1989. x+97 pp.

\bibitem{MY-10}
C.~G. Moreira and J.-C. Yoccoz.
\newblock Tangences homoclines stables pour des ensembles hyperboliques de
  grande dimension fractale.
\newblock {\em Annales Scientifiques de l'\'Ecole Normale Sup\'erieure},
  43:1--68, 2010, \url{http://www.numdam.org/articles/10.24033/asens.2115/}.

\bibitem{PT93}
J.~Palis and F.~Takens.
\newblock {\em Hyperbolicity and Sensitive chaotic dynamics at homoclinic
  biifurcations: fractal dimensios and infinitely many attractors}.
\newblock Cambridge Univ. Press, 1993.




\end{thebibliography}
\end{document}